\documentclass[11pt, a4paper]{amsart}
\usepackage{a4, a4wide}
\usepackage{amssymb}
\usepackage{amsmath}
\usepackage{amsthm}
\usepackage{amstext}
\usepackage{amscd}
\usepackage{enumitem}
\usepackage{latexsym}
\usepackage{graphics}
\usepackage{color}
\usepackage[all]{xy}
\usepackage{mathtools}

\usepackage[colorlinks, pagebackref]{hyperref}
\hypersetup{
  colorlinks=true,
  citecolor=blue,
  linkcolor=blue,
  urlcolor=blue}

\newtheorem{intro-thm}{Theorem}[]

\theoremstyle{plain}
\newtheorem{theorem}{Theorem}[section]

\newtheorem{lemma}[theorem]{Lemma}

\newtheorem{fact}[theorem]{Fact}
\theoremstyle{definition}
\newtheorem{remark}[theorem]{Remark}

\newtheorem{notation}[theorem]{Notation}

\newtheorem*{claim*}{Claim}

\numberwithin{equation}{section}
\newcommand{\ilim}{\mathop{\varprojlim}\limits} % inverse limit
\newcommand{\dlim}{\mathop{\varinjlim}\limits}  % direct limit

\def\~{\widetilde}
\def\-{\overline}
\def\<{\langle}
\def\>{\rangle}
\def\@{\mathcal}
\def\!{\mathbf}
\def\#{\mathbb}
\def\<{\langle}
\def\>{\rangle} 
\def\-{\overline} 
\def\~{\widetilde}
\def\^{\widehat}

\newcommand{\Proj}{{\rm Proj} \,}

\newcommand{\Spec}{{\rm Spec \,}}

\input{xy}
\xyoption{all}

\begin{document}

\title{Naive $\mathbb{A}^1$-connectedness of retract rational varieties}

\author{Chetan Balwe}
\address{Department of Mathematical Sciences, Indian Institute of Science Education and Research Mohali, Knowledge City, Sector-81, Mohali 140306, India.}
\email{cbalwe@iisermohali.ac.in}

\author{Bandna Rani}
\address{Department of Mathematical Sciences, Indian Institute of Science Education and Research Mohali, Knowledge City, Sector-81, Mohali 140306, India.}
\email{ph17049@iisermohali.ac.in}

%\subjclass[2010]{14C15, 14C25, 19E15 (Primary)}
%\keywords{}

\begin{abstract}
A smooth, proper, retract rational variety over a field $k$ is known to be $\#A^1$-connected. We improve on this result, in the case when $k$ is infinite, showing that such varieties are naively $\#A^1$-connected. 
\end{abstract}

\maketitle
%\tableofcontents

\section{Introduction}

Let $k$ be a field and let $Sm_k$ denote site of smooth varieties over $k$, equipped with the Nisnevich topology. To any simplicial sheaf $\@X$ on $Sm_k$, $\#A^1$-homotopy theory (see \cite{Morel-Voevodsky}) associates the sheaf of $\#A^1$-connected components of $\@X$, denoted by $\pi_0^{\#A^1}(\@X)$. We say that $\@X$ is $\#A^1$-connected if $\pi_0^{\#A^1}(\@X)$ is equal to $\ast$, the point sheaf. 

In general, the sheaf $\pi_0^{\#A^1}(\@X)$ is difficult to compute explicitly for general $\@X$. We recall a related construction --- the sheaf of naively $\#A^1$-connected components, which may be viewed as a crude approximation to $\pi_0^{\#A^1}(\@X)$ (see \cite{BHS} and \cite{BS1}). For any smooth scheme $U$, let $\sigma_0, \sigma_1: U \to U \times \#A^1$ denote the morphisms which map $U$ isomorphically onto the closed subschemes $U  \times \{0\}$ and $U \times \{1\}$ of $U \times \#A^1$, respectively.  If $\@X$ is a sheaf of sets, and $U$ is a smooth scheme over $k$, we say that two elements $f,g \in \@X(U)$ are $\#A^1$-homotopic if there exists an element $h \in \@X(U \times \#A^1)$ which is mapped to $f$ and $g$ under the two restriction maps $\sigma_0^*, \sigma_1^*: \@X(U \times \#A^1) \to \@X(U)$. This generates an equivalence relation on $\@X(U)$, which we denote by $\sim_U$. Let $\@S(\@X)$ be the Nisnevich sheafification of the presheaf which maps $U$ to the set of equivalence classes $\@X(U)/\sim_U$. We say that $\@X$ is naively $\#A^1$-connected if $\@S(\@X) = \ast$. 

The relationship between the functors $\pi_0^{\#A^1}(-)$ and $\@S(-)$ is quite subtle.  We have canonical morphisms of sheaves $\@X \to \@S(\@X)$ and $\@S(\@X) \to \pi_0^{\#A^1}(X)$, both of which are epimorphisms. It was proved in \cite{BHS} that there exists a canonical morphism 
\[
\pi_0^{\#A^1}(\@X) \to \dlim \@S^n(\@X) =: \@L(\@X)
\]
which is an isomorphism if and only if $\pi_0^{\#A^1}(X)$ is $\#A^1$-invariant. Thus, if $\@X$ is $\#A^1$-connected, then we have $\dlim \@S^n(\@X) = \ast$. A recent example, constructed by J. Ayoub shows that $\pi_0^{\#A^1}(\@X)$ need not be $\#A^1$-invariant in general (see \cite{Ayoub}). However, the map $\pi_0^{\#A^1}(\@X)(K) \to \@L(\@X)(K)$ is a bijection for any finitely generated, separable field extension $K/k$ (see \cite{BRS}. 

The relationship between $\pi_0^{\#A^1}(-)$ and $\@S$ is more interesting when we study it in the context of proper varieties. A result of Asok and Morel states that if $X$ is a proper variety over $k$, the map $\@S(X)(L) \to \pi_0^{\#A^1}(X)(L)$ is an isomorphism for every finitely generated, separable field extension $L/k$. A result of Morel shows that if $\pi_0^{\#A^1}(X)(L) = \ast$ for all finitely generated, separable field extensions $L/k$, then $X$ is $\#A^1$-connected. Thus, we see that if $X$ is a proper variety over $k$, then $X$ is $\#A^1$-connected if and only if $\@S(X)(L) = \ast$ for any finitely generated, separable field extension $L/k$. Since $X$ is proper over $k$, the latter condition is equivalent to saying that $X$ is universally $R$-trivial. 

The relationship between universal $R$-triviality and other near-rationality properties (such as retract rationality, unirationality, etc.) of smooth, proper varieties is not well-understood. (Note that smoothness is an important condition in this context. For example, a cone over any variety is $\#A^1$-connected, but need not have any other near-rationality property.) Smooth, proper, retract rational varieties over a field $k$ are $\#A^1$-connected. This was proved by Asok and Morel in the case $char(k) = 0$ (see \cite[Theorem 2.3.6]{Asok-Morel}) and in the general case by Kahn and Sujatha (see \cite[Theorems 8.5.1 and 8.6.2]{Kahn-Sujatha}). It is not known if the converse is true.  

We will prove the following result:
\begin{theorem}
\label{main theorem}
Let $k$ be an infinite field. Let $X$ be a smooth, proper, retract rational variety over $k$. Then $\@S(X) = \ast$. 
\end{theorem}

Thus, if $X$ is a smooth, proper variety over an infinite field $k$, we have the following implications:
\[
\text{retract rational} \implies \text{naively $\#A^1$-connected} \implies \text{$\#A^1$-connected}
\]
It is not known if either of these two implications is strict. We note that if $X$ is a proper, $\#A^1$-connected variety, a result of Sawant (see \cite{Sawant-IC}) implies that $\@S^2(X) = \ast$. It is easy to construct examples of singular varieties $X$ such that $\@S^2(X) = \ast$ but $\@S(X) \neq \ast$, but we are not aware of an example of a smooth variety satisfying these conditions. 

The above theorem says that if $X$ is a smooth, proper, retract rational variety over an infinite field $k$, then $\@S(X) \to \pi_0^{\#A^1}(X)$ is an isomorphism. This may be contrasted with the results of \cite{BS2} where we see that if $X$ is a smooth, proper surface that is birationally ruled over a curve of genus $>0$, the morphism $\@S(X) \to \pi_0^{\#A^1}(X)$ is an isomorphism if and only if $X$ is minimal. 

In Section \ref{section preliminaries}, we collect  some technical results regarding henselization that are required for the proof of Theorem \ref{main theorem}. At the beginning of the section, we include a brief, informal discussion regarding the strategy of the proof. The proof of Theorem \ref{main theorem} is presented in Section 3. 

\begin{notation}
\label{notation} \mbox{}
\begin{enumerate}
\item If $x$ is a point on a scheme $X$, the residue field at $x$ will be denoted by $\kappa(x)$. 

\item Given a scheme $X$ (resp. an affine scheme $X = \Spec R$), and an ideal sheaf $\@I$ (resp. an ideal $I \subset R$), we will denote by $\@Z(\@I)$ (resp. $\@Z(I)$) the closed subscheme of $X$ associated to the ideal sheaf $\@I$ (resp. the ideal $I$). If $\@L$ is a line bundle on $X$ and $S$ is a set of sections of $\@L$, we may also write $\@Z(S)$ for the closed subschemes defined by the vanishing of the elements of $S$. 

\item Let $\@F$ be a presheaf on $Sm_k$. Suppose $U$ is an essentially smooth scheme over $k$, i.e. a filtered inverse limit of a diagram of smooth schemes $U_{\alpha}$, in which the transition maps are \'etale, affine morphisms. Then, we will define $\@F(U) = \dlim \@F(U_{\alpha})$. We will primarily use this definition in the context of schemes of the form $\Spec \@O_{X,x}^h$ where $X$ is a smooth variety over $k$, $x$ is a point of $X$ and $\@O_{X,x}^h$ is the henselization of the local ring $\@O_{X,x}$. We will denote the local scheme $\Spec \@O_{X,x}^h$ by $X_x$ and the canonical morphism $X_x \to X$ by $\omega_x$.   
\end{enumerate}
\end{notation}

\section{Preliminaries}
\label{section preliminaries}

In this section, we gather some technical results required in our proof of Theorem \ref{main theorem}. First, we briefly discuss the role they will play in the proof. 

Let $k$ be a field and let $X$ be a variety over $k$. To prove that $\@S(X) = \ast$, it suffices to show that if $U$ is any smooth variety over $k$ and $u$ is any point of $U$, then $\@S(X)(U_u) = \ast$ (see Notation \ref{notation}). 

If $f: U_u \to X$ is any morphism, and if $f(u) = x$, then $f$ factors through the canonical morphism $\omega_x: X_x \to X$. (This follows from the definition of henselization in subsection \ref{subsection henselization}.) Thus, if we want to show that $\@S(X) = \ast$, it suffices to prove that $\@S(X)(X_x) = \ast$ for any point $x \in X$. In other words, we need to show that there exists a point $x_0 \in X(k)$ such that for any $x \in X$, there exists a chain of $\#A^1$-homotopies connecting $\omega_x$ to the morphism 
 \[
 X_x \to \Spec k \xrightarrow{x_0} X \text{.}
 \]

Thus, for any $x \in X$, we will need to construct certain $\#A^1$-homotopies of $X_x$ in $X$, i.e. morphisms of the form $\Spec R[t] \to X$ where $R = \@O_{X,x}^h$. For a general variety $X$, it may not be easy to construct such morphisms. However, we can start with a somewhat weaker notion. A morphism $\Spec R[[t]] \to X$ could be seen as the ``germ" of a homotopy (or an `infinitesimal homotopy", as in subsection \ref{subsection infinitesimal homotopies}). If $X = \#A^N_k$ for some $N$, there is a simple trick to obtain a homotopy from a morphism $\Spec R[[t]] \to X$. Indeed, this morphism is given by an $N$-tuple of power series with coefficients in $R$. We can just truncate these power series at some degree to obtain an $N$-tuple of polynomials. This gives us a morphism from $\Spec R[t]$ into $\#A^N_k$, i.e. an $\#A^1$-homotopy of $\Spec R$ in $\#A^N_k$. By truncating at a sufficiently high degree, we can make sure that the $\#A^1$-homotopy approximates the infinitesimal homotopy that we started with. With a little care, this trick can be made to work in $\#P^N_k$. 

Given a retract rational variety $X$, we have birational maps $X \dashrightarrow \#P^N_k$ and $\#P^N_k \dashrightarrow X$, the composition of which is the identity map on $X$. To prove that $X$ is naively $\#A^1$-connected, we have to construct some $\#A^1$-homotopies on $X$. In very crude terms, the idea is to first construct an infinitesimal homotopy on $X$ and then to ``lift it" to $\#P^N_k$, and to then use the ``truncation trick" described above to obtain an $\#A^1$-homotopy. 

Of course, we cannot really lift the infinitesimal homotopy to $\#P^N_k$, since we only have a rational map from $X$ to $\#P^N_k$. So one has to use a somewhat sophisticated analogue of the ``truncation trick" --- we need to divide the power series in question by an appropriately chosen Weierstrass polynomial and compute the remainders. For this, we will need an analogue of the Weierstrass preparation theorem in the context of henselian local rings. This is what we recall in subsection \ref{subsection henselization}. In subsection \ref{subsection infinitesimal homotopies} we prove a lemma that will allow us construct suitable infinitesimal homotopies on a smooth variety. This is precisely where smoothness is crucially required. Note that Theorem \ref{main theorem} is not true for singular varieties. (Indeed, a rational, proper variety that is not smooth is not necessarily $\#A^1$-connected.)

\subsection{Henselization}

\label{subsection henselization}

In this subsection, we will collect some facts about henselization that we will require in the proof of the main result. 

Let $R$ be a ring and let $I$ be an ideal of $R$. We say that $(R,I)$ is a \emph{henselian pair} if $I$ is contained in the Jacobson radical of $R$ and given any \'etale extension $R \to R'$ such that $R/I \to R'/IR'$ is an isomorphism, there exists an $R$-algebra homomorphism $R' \to R$. When $R$ is a local ring with maximal ideal $\mathfrak{m}$ such that $(R, \mathfrak{m})$ is a henselian pair, we say that $R$ is a \emph{henselian local ring}.

\begin{fact}[See {\cite[Theorem 18.5.1]{EGA4-4}}]
\label{fact henselian finite morphism}
Let $U = \Spec R$ where $R$ is a henselian local ring. Let $u$ be the closed point of $U$. Let $f: V \to U$ be a finite morphism. Then, $V$ is the disjoint union of the schemes $\Spec \@O_{V,v}$ where $v$ ranges over all the points of $f^{-1}(u)$. 
\end{fact}

Given any ring $R$ and an ideal $I$ of $R$, consider the category of henselian pairs $(S,J)$ such that $S$ is an $R$-algebra and $IS \subset J$. This category has an initial object $(R^h, I^h)$, called the \emph{henselization of $R$ at $I$}. When $R$ is a local ring with maximal ideal $\mathfrak{m}$, we will refer to the henselization of $R$ at $\mathfrak{m}$ as just the \emph{henselization of $R$}. 

\begin{fact}[See {\cite[Lemma 15.12.4]{StacksProject}}] \label{fact henselisation completion} Let $R$ be a noetherian ring and let $I$ be an ideal of $R$. Let $(R^h, I^h)$ denote the henselization of $R$ at $I$. Then the canonical homomorphism $R \to R^h$ induces an isomorphism of the $I$-adic completions $\^R \to \^{R^h}$. The canonical homomorphism $R \to R^h$ is a flat extension which induces an isomorphism of the $I$-adic completions $\^R \to \^{R^h}$. Also, $R^h$ is a noetherian ring and so $R^h \to \^R$ is a faithfully flat extension. 
\end{fact} 

We will primarily be concerned with the henselizations of the local rings at points of a variety over a field and we quote a basic result about such rings. 

\begin{fact}[See {\cite[Corollary 18.7.6]{EGA4-4}}]
\label{fact henselian local rings varieties}
Let $k$ be a field. Let $X$ be a variety over $k$ and let $x$ be a point of $X$. Then, the ring $\@O_{X,x}^h$ is an excellent ring.
\end{fact} 

Let $R$ be any ring and let $t$ denote the variable. The henselization of the ring $R[t]$ at the ideal $\<t\>$ is called the ring of henselian power series, and is denoted by $R\{t\}$. We will now review some results about such rings. 

First, let us take $R$ to be a noetherian ring. Then, Fact \ref{fact henselisation completion} implies that $R\{t\}$ injects into $R[[t]]$ and we will identify $R\{t\}$ with its image in $R[[t]]$. Since $R\{t\} \to R[[t]]$ is faithfully flat, it is easy to see that an element $f \in R\{t\}$ is a unit if and only if it is a unit in $R[[t]]$, i.e. if and only if its image under the quotient map $R[[t]] \to R[[t]]/\<t\> \cong R$ is a unit in $R$. 

The results stated in the previous paragraph continue to hold even when $R$ is not noetherian. Indeed, the functor $R \to R\{t\}$ commutes with filtered colimits. So one can express $R$ as a filtered colimit of its finitely generated $\#Z$-algebras and generalize the above statements to the non-noetherian setting. ( See \cite[Subsection 2.1.2]{Bouthier-Cesnavicius}.) 

If $(R, \mathfrak{m})$ is a local ring, a polynomial $f(t) \in R[t]$ is said to be a \emph{Weierstarss polynomial} if it is of the form $f(t) = t^d + a_{d-1}t^{d-1} + \ldots + a_0$ where $a_i \in \mathfrak{m}$ for all $i$. 

\begin{fact}
\label{fact Weierstrass preparation}
Let $(R, \mathfrak{m})$ be a henselian local ring. Let $f \in R\{t\} \backslash \mathfrak{m}R\{t\}$. Then $f$ can be uniquely factored as $f = P\cdot u$ where $f$ is a Weierstarss polynomial and $u$ is a unit  in $R\{t\}$. Also, the natural homomorphism $R[t]/\<P\> \to R\{t\}/\<f\>$ is an isomorphism. 
\end{fact}

The following lemma is an easy consequence: 

\begin{lemma}
\label{lemma preparation for ideals}
Let $R$ be a henselian local ring with maximal ideal $\mathfrak{m}$. Let $I$ be a proper ideal of $R[t]$ such that the following conditions hold:
\begin{enumerate}
\item The homomorphism $R \to R[t]/I$ is a finite extension. 
\item The only prime ideal of $R[t]$ containing $I$ and $\mathfrak{m}R[t]$ is $\<\mathfrak{m}, t\>$. 
\end{enumerate}
Then, the homomorphism $R[t]/I \to R\{t\}/IR\{t\}$ is an isomorphism.  
\end{lemma}
\begin{proof}
Let $Z = \Spec R[t]/I$, which we view as a closed subscheme of $\Spec R[t]$. Let $x_0$ be the closed point of $\Spec R$ and let $y_0$ be the point of $\Spec R[t]$ corresponding to the ideal $\<\mathfrak{m}, t\>$. 

Let $\pi: Z \to \Spec R$ be the morphism corresponding to the $R$-algebra homomorphism $R \to R[t]/I$. According to condition (1), $\pi$ is a finite morphism. Thus, if $z \in Z$ is any point, there exists a point $z_0$ in its closure such that $\pi(z_0) = x_0$. By (2), the closed subscheme $Z$ and the fibre $\pi^{-1}(x_0)$ have only the point $y_0$ in common. Thus, we see that every point of $Z$ lies in the closure of $y_0$. 

By Fact \ref{fact henselian finite morphism}, we see that $Z$ is isomorphic to $\Spec \@O_{Z,y_0}$. Thus, we see that if $S = R[t]_{\<\mathfrak{m}, t\>}$, then the homomorphism $R[t]/I \to S/IS$ is an isomorphism. 

The ring $R\{t\}$ is a local ring with maximal ideal $\<\mathfrak{m}, t\>$. Thus, the canonical homomorphism $R[t] \to R\{t\}$ induces a local homomorphism $S \to R\{t\}$. This homomorphism is flat, and since it is a local homomorphism, it is faithfully flat. Thus, as $IS \not \subset \mathfrak{m}S$, we see that $I R \{t\} \not\subset \mathfrak{m} R\{t\}$. Let $f$ be an element of $I R\{t\} \backslash \mathfrak{m} R \{t\}$. Then $f = u \cdot p$ where $u$ is a unit in $R\{t\}$ and $p$ is a Weierstrass polynomial. Since $S \to R\{t\}$ is a faithfully flat extension, $IS = S \cap IR\{t\}$. Thus, $p \in IS$. 

By \cite[Proposition 3.1.2]{Bouthier-Cesnavicius}, the ring homomorphism $R[t]/pR[t] \to R\{t\}/pR\{t\}$ is an isomorphism. Thus, it follows that the homomorphism $S/pS \to R\{t\}/pR\{t\}$ is surjective. It is also injective since it is a faithfully flat extension. Thus, it is an isomorphism. 

As $pS \subset IS$ and $pR\{t\} \subset I R\{t\}$, it follows that $S/IS \to R\{t\}/IR\{t\}$ is an isomorphism. This completes the proof. 
\end{proof}

A noetherian local ring $(R, \mathfrak{m})$ is said to be an \emph{approximation ring} if for any finite system of polynomial equations with coefficients in $R$, the set of solutions in $R$ is dense, with respect to the $\mathfrak{m}$-adic topology, in the set of solutions in the $\mathfrak{m}$-adic completion $\widehat{R}$. We will need the following facts: 

\begin{fact}[See {\cite[Theorem 1.3]{Popescu}}] Excellent henselian local rings are approximation rings. 
\label{fact popescu 1}
\end{fact}

\begin{fact}[See {[\cite[Corollary 3.5]{Popescu}}] Let $(R, \mathfrak{m})$ be an approximation ring. Then $R\{t\}$ is an approximation ring. 
\label{fact popescu 2}
\end{fact}

\subsection{Infinitesimal homotopies}

\label{subsection infinitesimal homotopies}

For any ring $R$, an \emph{infinitesimal homotopy} of $\Spec R$ in a scheme $X$ is a morphism $h: \Spec R\{t\} \to X$. Let $\^\sigma_0: \Spec R \to \Spec R\{t\}$ be the morphism induced by the quotient homomorphism $R\{t\} \to R\{t\}/tR\{t\} \cong R$.
We say that this infinitesimal homotopy \emph{starts} from the morphism $h \circ \^\sigma_0: \Spec R \to X$. 

The following lemma shows that for a smooth variety $X$ and a point $x \in X$, one can easily construct the germ of a homotopy of $X_x$ in $X$ starting from the canonical morphism $\omega_x: X_x \to X$. 

\begin{lemma}
\label{lemma infinitesimal homotopy}
Let $X$ be a smooth $d$-dimensional variety over $k$. Let $x$ be a point of $X$. Let $U$ be an open subset of $X$. Let $R = \@O_{X,x}^h$ and let $\omega_x: \Spec(R) =: X_x \to X$ be the canonical morphism. Then, there exists a morphism $h: \Spec R\{t\} \to X$ starting from $\omega_x$ such that 
\[
\Spec \kappa(x)\{t\} \to \Spec R\{t\} \stackrel{h}{\to} X
\]
maps the generic point of $\Spec \kappa(x)\{t\}$ into $U$. 
\end{lemma}

\begin{proof}
If $x \in U$, we may take $h$ to be the composition 
\[
\Spec R\{t\} \to \Spec R \stackrel{\omega_x}{\to} X
\]
where the first morphism is induced by the inclusion $R \hookrightarrow R\{t\}$. Thus, we will now assume that $x \notin U$. Let $Z = X \backslash U$. 

For any non-negative integer $n$, consider the functor 
\[
U \mapsto Mor_{Sch/k}(U \times \Spec k[t]/\<t^{n+1}\>) 
\]
on the category of $k$-schemes. This functor is known to be representable by a $k$-scheme of finite type (see \cite{Greenberg}), which we denote by $J_n(X)$. Let $J(X) = \ilim J_n(X)$, where the inverse limit is computed in the category of $k$-schemes. The quotient homomorphism $k[[t]] \to k[t]/\<t^{n+1}\>$ induces the morphism $\pi_{n}^X: J(X) \to J_n(X)$. For $n \geq m$, let $\pi^X_{n,m}: J_n(X) \to J_m(X)$ denotes the morphism induced by the quotient homomorphism $k[t]/\<t^{n+1}\> \to k[t]/\<t^{m+1}\>$. 

Choose a morphism $\gamma: \Spec \kappa(x)[[t]] \to X$ which maps the closed point of $\Spec \kappa(x)[[t]]$ to $x$ and the generic point into $U$. For any $n \geq 0$, let $\gamma_n$ denote the composition 
\[
\Spec \kappa(x)[t]/\<t^{n+1}\> \to \Spec \kappa(x)[[t]] \xrightarrow{\gamma} X \text{.}
\]
We identify $\gamma_n$ with a $\kappa(x)$-valued point of $J_n(X)$, which we denote by $\~{\gamma}_n$. We define $\~g_0: \Spec R \to J_0(X) = X$ to be the morphism $\omega_x$. For $n \geq 1$, we will inductively construct a morphism $\~g_n: \Spec R \to J_n(X)$ such that: 
\begin{itemize}
\item[(i)] the composition $\Spec \kappa(x) \to \Spec R \xrightarrow{\~g_n} J_n(X)$
is equal to $\~{\gamma}_n$, and 
\item[(ii)] the composition $\pi^X_{n+1,n} \circ \~g_{n+1}$ equals $\~g_n$. 
\end{itemize} 
Suppose $\~g_n$ has been chosen for some non-negative integer $n$. Since $X$ is smooth, $J_{n+1}(X) \to J_n(X)$ is smooth. (In fact, it is an affine bundle for all $n$ --- see \cite[Lemma 9.1]{Looijenga}). So we can choose a morphism $\~g_{n+1}: \Spec R \to J_{n+1}(X)$ satisfying the conditions (i) and (ii) (see \cite[Corollary 17.16.3, (ii)]{EGA4-4}). 

The collection $\{\~g_n\}_{n \geq 0}$ defines a morphism $\~g: \Spec R \to J(X)$, which corresponds to a morphism $g: \Spec R[[t]] \to X$. The restriction of $g$ to $\kappa(x)[[t]]$ is equal to $\gamma$. 

There exists an integer $n$ such that if $\gamma': \Spec \kappa(x)[[t]] \to X$ satisfies $\pi_{n, \Spec \kappa(x)}(\gamma) = \pi_{n,\Spec \kappa(x)}(\gamma')$, then $\gamma' $ maps the generic point of $\Spec \kappa(x)[[t]]$ into $U$. (Indeed, if there is no such $n$, then since $J(Z) = \ilim J_n(X)$, it will follow that $\gamma \in J(Z)$, which is not true.) By Facts \ref{fact popescu 1} and \ref{fact popescu 2}, there exists a morphism $h: \Spec R\{t\} \to X$ such that $\pi^X_n(h) = \pi^X_{n}(g)$. This proves the lemma. 
\end{proof}

\begin{remark}
As we see in the above proof, it is very easy to construct a morphism $\Spec R[[t]] \to X$. We could have used this as our notion of ``infinitesimal homotopy", if we had an analogue of Fact \ref{fact Weierstrass preparation} for the ring $R[[t]]$, at least when $R$ is regular. (Such a result was proved in characteristic $0$ by Lafon in \cite{Lafon}.) The proof of the preparation theorem for $R\{t\}$ in \cite{Bouthier-Cesnavicius} crucially uses the fact that the functor $R\{t\}$ is a colimit of finite type $R$-algebras and we do not know if it can be adapted to give an analogous result for $R[[t]]$. So we choose to work with the ring $R\{t\}$ instead. 
\end{remark}

\section{Retract rational varieties}

We will now focus on retract rational varieties. In subsection \ref{subsection rational curves}, we prove a small lemma regarding rational curves in projective space. The proof of Theorem \ref{main theorem} is presented in subsection 3.2. 

\subsection{Rational curves in projective space}
\label{subsection rational curves}

Let us fix a base field $k$. Let $L$ be any field containing $k$. We will use $T_0$ and $T_1$ as homogeneous coordinates on $\#P^1_L$. In other words, we will write $\#P^1_L = \Proj L[T_0, T_1]$.  We will identify $\#A^1_L = \Spec L[t]$ with the open subscheme $\#P^1_L \backslash \@Z(T_1)$ by identifying $t$ with $T_0/T_1$. We will denote the point $(0:1)$ of $\#P^1_L$ by $\mathbf{0}_L$ and the point $(1:0)$ by $\mathbf{\infty}_L$. Thus, $t$ is a parameter at $\mathbf{0}_L$ and $1/t$ is a parameter at $\mathbf{\infty}_L$. To avoid making the notation cumbersome, we will write $\mathbf{0}$ and $\mathbf{\infty}$ instead of $\mathbf{0}_L$ and $\mathbf{\infty}_L$ in the following discussion. 

A morphism $\phi: \#P^1_L \to \#P^N_k$ can be represented by an $(N+1)$-tuple $(P_0, \ldots, P_N)$ of homogeneous polynomials of a fixed degree $d$, such that $P_i \neq 0$ for some $i$. (Of course, some of the $P_i$'s may be equal to $0$. The zero polynomial can be assigned any degree.) Such a representation is not unique, but if we require the polynomials to be coprime, it is unique up to multiplication by a unit. Dehomogenizing this $(N+1)$-tuple with respect to $T_1$ gives an $(N+1)$-tuple of polynomials in $t$, which describes the restriction of $\phi$ to the open subscheme $D(T_1)$. 

Recall that given a morphism from $\#P^1_L \backslash \@Z(T_1)$ to $\#P_k^N$, it can be uniquely extended to a morphism $\#P^1_L \to \#P^N_k$. A morphism from $\#P^1_L \backslash \@Z(T_1)$ to $\#P^N_k$ is given by an $(N+1)$-tuple of polynomials in $L[t]$, such that at least one of the polynomials is non-zero. Thus, we see that a morphism $\#P^1_k \to \#P^N_k$ can be represented in three ways --- using $(N+1)$-tuples of homogeneous polynomials of a same degree in $(T_0,T_1)$ or by using $(N+1)$-tuples of polynomials in either $t$ or $1/t$. (Again, note that these representations are unique up to multiplication by a unit if we require the polynomials to be coprime.)

Given a morphism $\phi: \#P^1_k \to \#P^N_k$, we choose a representation of $\phi$ by an $(N+1)$-tuple of polynomials $(P_0, \ldots, P_N)$ in $k[t]$ which are coprime. Let $m \geq 0$ be any integer. For $0 \leq i \leq N$, let $P'_i$ be the polynomial obtained by truncating $P_i$ to degree $m$. Then, the $(N+1)$-tuple $(P'_0, \ldots, P'_N)$ is called the $m$-jet of $\phi$ at $\mathbf{0}$. Similarly, we can define the $m$-jet of $\phi$ at $\mathbf{\infty}$. Note that these are well-defined up to multiplication by a unit. 

The following lemma shows that given a closed subscheme $W$ of $\#P_k^N$ of codimension $\geq 2$, and non-negative integers $m_1$ and $m_2$, there exists a morphism $\phi: \#P^1_k \to \#P^N_k$ such that it maps $\#P^1_k \backslash \{\mathbf{0}, \mathbf{\infty}\}$ into $\#P^N_k \backslash W$, has a prescribed $m_1$-jet at $\mathbf{0}$ and a prescribed $m_2$-jet at $\mathbf{\infty}$. 

\begin{lemma}
\label{lemma constructing rational curve}
Let $k$ be an infinite field and let $L$ be a field containing $k$. Let $N$ be a positive integer. Let $W$ be a closed subscheme of $\#P^N_k$ of codimension $\geq 2$. Let $P = (P_0, \ldots, P_N)$ and $(Q_0, \ldots, Q_N)$ be $(N+1)$-tuples of polynomials in $L[t]$. Assume that $P_i \neq 0$ and $Q_j \neq 0$ for some indices $i,j$. Let $m_1$ be an integer such that $m_1 \geq \max_i deg P_i$. Assume that $P_i(0) \neq 0$ for some $i$. Then, there exists an $(N+1)$-tuple $(c_0, \ldots, c_N) \in k^{N+1}$  such that the following conditions hold:
\begin{itemize}
\item[(a)] For $0 \leq i \leq N$, let $R_i(t)= P_i(t) + t^{m_1+1}c_i + t^{m_1 + 2} Q_i(t)$. Then, the polynomials $R_0, \ldots, R_{N}$ are coprime. 
\item[(c)] Let $\phi: \#P^1_L \to \#P^N_k$ be the morphism represented by the $(N+1)$-tuple $(R_0, \ldots, R_N)$. Then $\phi(\#P^1_L \backslash \{\mathbf{0}, \mathbf{\infty} \}) \subset \#P^N_k \backslash W$. 
\end{itemize}
\end{lemma}
Note that if $\phi$ is as described in the lemma, the $m_1$-jet of $\phi$ at $\mathbf{0}$ is $(P_0(t), \ldots, P_N(t))$ and if $m_2 = \max_i \deg Q_j$, then the $m_2$-jet of $\phi$ at $\infty$ is $(Q_0(1/t), \ldots Q_N(1/t))$. 

\begin{proof} For any point $x = (x_0, \ldots, x_N)$ of $\#A^{N+1}_k$ and $0 \leq i \leq N$, we define  
\[
R^x_i(t) = P_i(t) + x_i t^{m_1 + 1} +  Q_i(t) t^{m_1 + 2}\text{.}
\]
Let $\phi_x: \#A^1_{\kappa(x)} \to \#A^{N+1}_k$ be the morphism defined by
\[
\phi_x(s) = (R^x_0(s), \ldots, R^x_N(s)) 
\]
for any $s \in \#A^1_{\kappa(x)}$. Let $C(W) \subset \#A^{N+1}$ be the cone over $W$. 

We define
\[
B := \{(x,s,z) | \phi_x(s) = z\} \subset  \#A^{N+1}_k \times (\#A^1_k \backslash \{0\}) \times C(W)\text{.}
\]
This is a closed subset of $\#A^{N+1}_k \times (\#A^1 \backslash \{0\}) \times C(W)$. Let $pr_1: B \to \#A^N_k$ be the projection map onto the first factor. We need to show that the complement of the image of $pr_1$ contains some $k$-rational point. 

Let $pr_{23}: \#A^{N}_k \times (\#A^1_k \backslash \{0\}) \times C(W) \to (\#A^1_k \backslash \{0\}) \times C(W)$ by the projection map onto the product of the second and third factors. We would like to estimate the dimension of the fibre $pr_{23}^{-1}(\gamma)$ where  $\gamma = (s,z) \in \#A^1_k \backslash \{0\} \times C(W)$. The equation $\phi_x(s) = z$ imposes $N+1$ linear conditions on $\#A^{N+1}_k$. Thus, the fibre has dimension $0$. Thus, it follows that $\dim(B) \leq 0 + 1 + \dim(C(W)) \leq N$.

It follows that the closure of $pr_1(B)$ is of dimension $ \leq N$. The result follows since $k$ is an infinite field. 
\end{proof}

\begin{remark}
This lemma is one of the main reasons for requiring the field $k$ to be infinite in Theorem \ref{main theorem}. The lemma need not hold if $k$ is a finite field since all the $k$-rational points of $\#P^N_k$ may be contained in $W$. 
\end{remark}

\subsection{Naive $\#A^1$-connectedness of retract rational varieties}
\label{subsection proof of main theorem}

We will now prove the main result of this paper. 

\begin{theorem}
\label{theorem main theorem section 3}
Let $k$ be an infinite field. Let $X$ be a smooth, proper, retract rational variety over $k$. Then $\@S(X) = \ast$. 
\end{theorem}
\begin{proof} 
Since $X$ is retract rational, there exists a positive integer $N\geq 1$, and rational maps $\phi: X \dashrightarrow \#P_k^N$ and $\psi: \#P_k^N \dashrightarrow X$ such that $\psi \circ \phi$ is the identity map on $X$. Since $k$ is infinite, this implies that $X(k)$ is non-empty. Since $X$ is a smooth, proper, retract rational variety, we have $\pi_0^{\#A^1}(X) = \ast$, and so $\@S(X)(k) = \ast$. Thus, to prove that $X$ is naively $\#A^1$-connected, it suffices to prove that for any point $x$, there exists a chain of $\#A^1$-homotopies of $X_x$ in $X$ connecting the canonical morphism $\omega_x: X_x \to X$ to a morphism that factors through some morphism $\Spec k \to X$. 

Let us fix a point $x \in X$. We will denote the ring $\@O_{X,x}^h$ by $R$. As in Notation \ref{notation}, we denote $\Spec R$ by $X_x$ and the canonical morphism $X_x \to X$ by $\omega_x$. We now set up some notation for working with the scheme $\#P^1_R$. 

We use the notation in subsection \ref{subsection rational curves}, so that $\mathbf{0}$ and $\mathbf{\infty}$ denote the points $(0:1)$ and $(1:0)$ of $\#P^1_k = \Proj k[T_0, T_1]$ respectively. Let $\sigma_0$ and $\sigma_{\infty}$ be the sections of the projection morphism $\#P^1_{R} \cong \#P^1_k \times \Spec R \to \Spec R$, mapping $\Spec R$ isomorphically onto the closed subschemes $\@Z(T_0) = \{\mathbf{0}\} \times \Spec R$ and $\@Z(T_1) = \{\mathbf{\infty}\} \times \Spec R$, respectively. We will denote the rational function $T_0/T_1$ by $t$ and thus identify the open subscheme $\#P^1_R \backslash \@Z(T_1)$ with $\Spec R[t]$. 

We will construct a morphism $H: \#P^1_R \to X$ such that $H \circ \sigma_0 = \omega_x$ and $H \circ \sigma_{\infty}$ factors through some morphism $\Spec k \to X$. Clearly, this will prove the result.  

There exists an ideal sheaf $\@K$ on $\#P^N_k$ such that if $\pi: Y \to X$ is the blowup of $\#P_k^N$ at $\@K$,  the map $\chi := \psi \circ \pi$ is a morphism from $Y$ to $X$. The sheaf $\@K$ can be chosen so that $W := \@Z(\@K)$ is a variety of codimension $\geq 2$. Let $V = \#P_k^N \backslash W$. Let $U \subset X$ be an open subset on which $\phi$ is defined and such that $\phi(U) \subset V$. The ideal sheaf $\@K$ corresponds to a homogeneous ideal of $k[X_0, \ldots, X_N]$ generated by homogeneous polynomials $p_1, \ldots, p_r$. We may assume, without loss of generality that the polynomials $p_1, \ldots, p_r$ are all of the same degree. Note that $r \geq 1$. 

The polynomials $p_1, \ldots, p_r$ define global sections of $\@K$, which generate $\@K$. Thus, we obtain a surjective morphism $\@O_{\#P^N}^r \to \@K$. Thus, we have the following sequence of homomorphisms of sheaves of graded $\@O_{\#P^N}$-rings
\[
Sym(\@O^r_{\#P^N}) \to Sym(\@K) \to \bigoplus_{j=0}^{\infty} \@K^j \text{.}
\]
This gives us a closed embedding of $Y$ into $\#P_k^N \times \#P_k^{r-1}$. 

Using Lemma \ref{lemma infinitesimal homotopy}, we choose an infinitesimal homotopy $h: \Spec R\{t\} \to X$ starting at $\omega_x$ such that $h$ maps the point $\eta_0$ of $\Spec R\{t\}$, corresponding to the ideal $\mathfrak{m}R\{t\}$, into $U$. This gives us a rational map $\phi \circ h: \Spec R \{t\} \dashrightarrow \#P^N$, which can be represented by an $N$-tuple $\!{f}:= (f_0, \ldots, f_N)$ where $f_i \in R\{t\}$. We choose the $f_i$ to be coprime in the unique factorization domain $R\{t\}$. Let $I$ denote the ideal $\<f_0, \ldots, f_N\>$. The rational map $\phi \circ h$ is a morphism if and only if this ideal is principal. Let  $J$ denote the ideal $\<p_1(\!f), \ldots, p_r(\!f)\>$.

Recall that we have chosen $h$ in such a way that the point  $\eta_0$ of $\Spec R\{t\}$, corresponding to the ideal  $\mathfrak{m} R \{t\}$, is mapped into $U$. Thus, the rational map $f \circ h$ is well-defined on $\eta_0$. Since we have chosen the $f_i$ to be coprime elements of the unique factorization domain $R\{t\}$, it follows that at least one of the $f_i$ does not vanish on $\eta_0$. By performing a change of coordinates on $\#P^N$, we may reduce to the situation where none of the $f_i$ vanishes on $\eta_0$. (Such a change of coordinates exists since $k$ is an infinite field.) Thus, for each $i$, we have $f_i = u_i \~f_i$ where $u_i$ is a unit in $R \{t\}$ and $\~f_i$ is a Weierstrass polynomial. 

Recall that $f$ maps the point $\eta_0$ into $V$. The zero set of the ideal $\<p_1, \ldots, p_r\>$ is contained in the complement of $V$. Thus, at least one of the polynomials $p_i$ does not vanish on $\eta_0$. We may assume that all the $p_i$'s are of the same degree. Suppose $p_1$ does not vanish on $\eta_0$. Then, for each $i \neq 1$, we can replace $p_i$ by $p_i  + \epsilon_i p_1$ where $\epsilon_i$ is $0$ if $p_i$ does not vanish at $\eta_0$ and is equal to $1$ otherwise. Thus, we may assume that none of the polynomials $p_i$ vanishes at $\eta_0$. 

Thus, $p_i(\!f) = v_i \cdot P_i$ where $v_i$ is a unit in $R\{t\}$ and $P_i$ is a Weierstrass polynomial. Let $p := t \cdot \left( \prod_i \~f_i \right)^2 \cdot \left( \prod_j P_j \right)^2$. Then, $p$ is an element of the ideal $IJ$ and it is a Weierstrass polynomial with $\deg_t(p) \geq 1$. 

For $0 \leq i \leq N$, we can express $f_i$ in the form 
\[
f_i = \alpha_i + p \beta_i
\]
where $\alpha_i  \in R[t]$ with $\deg_t (\alpha_i) < \deg_t(p)$ and $\beta_i \in R\{t\}$. 

For $0 \leq i \leq N$, let $\-{\alpha}_i(t) \in \kappa(x)[t]$ be the image of $\alpha_i(t)$ under the quotient homomorphism $R[t] \to R[t]/\mathfrak{m}R[t] = \kappa(x)[t]$. Let $d$ be the largest non-negative integer such that $t^d$ divides $\-{\alpha}_i(t)$ for all $i$. Let $(\lambda_0, \ldots, \lambda_N) \in k^{N+1}$ such that $(\lambda_0: \ldots: \lambda_N) \in V(k)$. (The existence of such an $(N+1)$-tuple $(\lambda_0, \ldots, \lambda_N)$ follows from the assumption that $k$ is infinite.) We apply Lemma \ref{lemma constructing rational curve} to the two $(N+1)$-tuples 
\[
(\-{\alpha}_0(t)/t^d, \-{\alpha}_1(t)/t^d, \ldots, \-{\alpha}_N(t)/t^d) \hspace{1cm} \text{and} \hspace{1cm} (\lambda_0, \ldots, \lambda_N)
\]
of polynomials in $\kappa(x)[t]$. We see that there exists an $(N+1)$-tuple $(\mu_0, \ldots, \mu_N) \in k^{N+1}$ such that the $(N+1)$ polynomials $R_0(t), \ldots, R_N(t)$ in $\kappa(x[t]$ the polynomials defined by
\[
R_i(t) = \-{\alpha}_i(t)/t^d + (\mu_i + \lambda_i t) \cdot t^{\deg_t(p)-d} \text{,}
\]
are coprime and define a morphism $\phi: \#P^1_{\kappa(x)} \to \#P^N_k$ mapping $\#P^1_{\kappa(x)} \backslash \{\mathbf{0}, \mathbf{\infty}\}$ into $\#P^N_k \backslash W$. 

Let $\!g = (g_0, \ldots, g_N)$ where $g_i(t) = \alpha_i(t) + (\mu_i + \lambda_i t) \cdot p(t)$ for $0 \leq i \leq N$. Then, if $\-g_i(t) \in \kappa(x)[t]$ is the image of $g_i(t)$ under the quotient homomorphism $R[t] \to \kappa(x)[t]$, we see that $\-g_i(t) = t^d R_i(t)$ for all $i$. Thus, the following conditions hold: 
\begin{itemize}
\item[(A)] The polynomials $\-g_0(t), \ldots, \-g_N(t)$ have no common zeros in $\#A^1_{\kappa(x)} \backslash \{0\}$.
\item[(B)] If $\-{\!g}$ denotes the $(N+1)$-tuple $(\-g_0, \ldots, \-g_N)$, then the polynomials $p_1(\-{\!g}), \ldots, p_r(\-{\!g})$ have no common zero in $\#A^1_{\kappa(x)} \backslash \{0\}$. 
\end{itemize}
Let $\~I$ and $\~J$ be the ideals $\<g_0, \ldots, g_N\>$ and $\<p_1(\!g), \ldots, p_r(\!g)\>$ of $R[t]$ respectively.

For $0 \leq i \leq N$,  
\begin{equation}
\label{equation: good approximation 1}
g_i = f_i + p \cdot (\mu_i + \lambda_i t - \beta_i) = f_i[1 + (p/f_i)\cdot (\mu_i + \lambda_it  - \beta_i)] \text{.} 
\end{equation}
As $t$ divides $p/f_i$, it is a non-unit in $R\{t\}$, and so $g_i$ is a unit multiple of $f_i$ in $R\{t\}$, for each $i$ . In particular, we have $\~I R\{t\} = I$.

Similarly, for $1 \leq i \leq r$, 
\[
p_i(\!g) - p_i(\!f) = \sum_j(g_j - f_j) \cdot Q_{ij}
\]
where $Q_{ij}$ is some element of $R\{t\}$. Thus,
\begin{equation}
\label{equation: good approximation 2}
p_i(\!g) = p_i(\!f)[1 + (p/p_i(\!f)) \cdot Q'_i] 
\end{equation}
for some $Q'_i \in R\{t\}$. As $t$ divides $p/p_i(\!f)$, it is a non-unit in $R\{t\}$,  and so $p_i(\!g)$ is a unit multiple of $p_i(\!f)$ in $R\{t\}$. This proves that $\~J R\{t\}= J$. 

We would now like to show that the $R$-algebra homomorphism $R[t]/(\~I\~J) \to R\{t\}/(IJ)$ is an isomorphism. This statement is trivially true if $\~I\~J$ is the unit ideals. Thus, let us assume for now, that $\~I\~J$ is not the unit ideal. We would like to apply Lemma \ref{lemma preparation for ideals}, and so we verify that $\~I\~J$ satisfies the hypothesis of that lemma. 

Condition (A), which was imposed on the $(N+1)$-tuple $\!g$, implies that if $\~I$ is not the unit ideal, then the only prime ideal containing $\~I$ and $\mathfrak{m}R[t]$ is $\<m,t\>$. Condition (B) implies the same for $\~J$. Since at least one of the ideals $\~I$ and $\~J$ is not the unit ideal, it follows that $\~I\~J$ satisfies condition (2) of Lemma \ref{lemma preparation for ideals}.

Now, we verify that the homomorphism $R \to R[t]/(\~I\~J)$ is a finite extension. For this, it suffices to find an element of $\~I \~J$ such that its leading coefficient, (as a polynomial in $t$) is an invertible element of $R$. 

First, we note that there exists an index $i_0$, $0 \leq i_0 \leq N$ such that $\lambda_{i_0} \neq 0$. Thus, $g_{i_0}$ is a polynomial in $t$ with a leading coefficient that is a unit in $R$. Secondly, we observe that for $1 \leq j \leq r$, $p_j(\!g)$ is a polynomial in $t$ with degree $\leq \deg(p_j) \cdot (\deg_t(p) + 1)$. The coefficient of $t^{\deg(p_j) \cdot (\deg_t(p) + 1)}$ is equal to $p_j(\lambda_0, \ldots, \lambda_N)$. Since $(\lambda_0: \ldots, \lambda_N)$ is in $V$, there exists an index $j_0$ such that $p_{j_0}(\lambda_0 , \ldots, \lambda_N) \neq 0$. Thus, $p_{j_0}(\!g)$ is a polynomial in $t$ with a leading coefficient that is an invertible element in $R$. Thus, $g_{i_0} \cdot p_{j_0}(\!g) \in \~I\~J$ is a polynomial in $t$ with a leading coefficient that is an invertible element of $R$. Thus, we see that the homomorphism $R \to R[t]/(\~I\~J)$ is a finite extension. 

Thus, if $\~I \~J$ is not the unit ideal, we may now apply Lemma \ref{lemma preparation for ideals} to conclude that the morphism 
\[\theta: \Spec R\{t\} \to \Spec R[t]
\]
induces an isomorphism of the closed subscheme $\@Z(IJ) \subset \Spec R\{t\}$ with the closed subscheme $\@Z(\~I \~J) \subset \Spec R[t]$. Of course, as we noted above, if $\~I\~J$ happens to be the unit ideal, $\@Z(IJ) \to \@Z(\~I\~J)$ is trivially an isomorphism. 

Let $\tau: B \to \Spec R\{t\}$ denote the blowup of $\Spec R\{t\}$ at the ideal $IJ$ and let $\~\tau: \~B \to \Spec R[t]$ denote the blowup of $\Spec R[t]$ at the ideal $\~I \~J$. Since $\~I\~JR\{t\} = IJ$, we see that 
\[
B \cong \~B \times_{\Spec R[t]} \Spec R \{t\} \text{.}
\] 
Let us denote the projection morphism $B \to \~B$ by $\theta'$. Since $\theta$ maps $\@Z(IJ)$ isomorphically onto the closed subscheme $\@Z(\~I\~J)$ of $\Spec R[t]$, it follows that $\theta'$ maps $\tau^{-1}(\@Z(IJ))$ isomorphically onto $\~\tau^{-1}(\@Z(\~I\~J))$. 

Since the ideal sheaf $\tau^{-1}(I) \cdot \@O_B$ is invertible, the rational map $h': \Spec R\{t\} \to \#P^N_k$ defined by the $(N+1)$-tuple $(f_0, \ldots, f_N)$ lifts to a morphism $B \to \#P^N$. Since the ideal sheaf $\tau^{-1}(J) \cdot \@O_B$ is also invertible, it further lifts to a morphism $h'': B \to Y$. Thus, the diagram 
\[
\xymatrix{
B \ar[r]^{h''} \ar[d]_{\tau} & Y \ar[d]^{\chi} \\
\Spec R\{t\} \ar[r]^h & X
}
\]
commutes. 

Similarly, the rational map $\~h': \Spec R[t] \dashrightarrow \#P^N_k$, defined by the $(N+1)$-tuple $(g_0, \ldots, g_N)$ lifts to a morphism $\~h'': \~B \to Y$. Thus, the diagram 
\[
\xymatrix{
B \ar[r]^{\theta'} \ar[d]_{\tau} & \~B \ar[r]^{\~h''} \ar[d]_{\~\tau} & Y \ar[d]^{\chi} \\
\Spec R\{t\} \ar[r]^{\theta} &\Spec R[t]  & X
}
\]
commutes. Notice that, in the above diagram, we do not yet have a morphism from $\Spec R[t]$ to $X$ making the diagram commute. We will prove that such a morphism exists. 

We have the two morphisms $h''$ and $\~h'' \circ \theta'$ from $B$ to $Y$. These need not be equal. However we will show that they agree on $\tau^{-1}(\@Z(IJ))$. This claim is trivial if $IJ$ is the unit ideal. Thus, we now assume that $IJ$ is not the unit ideal. 

Let $z$ be any point of $\tau^{-1}(\@Z(IJ))$. We want to prove that that $h''(z) = \~h'' \circ \theta'(z)$. Recall that we have fixed an embedding of $Y$ into $\#P_k^N \times \#P_k^{r-1}$. Let $pr_1: \#P_k^N \times \#P_k^{r-1} \to \#P_k^N$ and $pr_2: \#P_k^N \times \#P_k^{r-1} \to \#P_k^{r-1}$ be the projection morphisms. It will suffice to prove that $pr_i \circ h''(z) = pr_i \circ \~h'' \circ \theta'(z)$ for $i = 1,2$. 

For any element $r \in R\{t\}$, we will denote its image in $\@O_{B,z}$ by $r$ as well. Let $\mathfrak{n}_z$ denote the maximal ideal of the local ring $\@O_{B,z}$. Since $\tau(z) \in \@Z(IJ)$, there exists at least one element in the set $\{f_0, \ldots, f_N, p_1(\!f), \ldots, p_r(\!f)\}$ which is a non-unit in $\@O_{B,z}$. Let us pick one such element and denote it by $q_z$. Recall that we had chosen $p$ to be equal to $t \cdot \left( \prod_i \~f_i \right)^2 \cdot \left( \prod_j P_j \right)^2$. Thus, it follows that, in the ring $\@O_{B,z}$, the non-unit element $q_z$ divides $p/f_i$ for every $i$ and $p/p_j(\!f)$ for every $j$. We will use this observation in the following discussion. 

The restriction of $pr_1 \circ h''$ to $\Spec \@O_{B,z}$ is given by the $(N+1)$-tuple $(f_0, \ldots, f_N)$.  We know that the ideal $I \cdot \@O_{B,z}$ is principal. Thus,   there exists an index $i_z$ such that $f_i/f_{i_z} \in \@O_{B,z}$ for all $i$. Let $f'_i = f_i/f_{i_z}$ for $0 \leq i \leq N$. If $\-{f'_i}$ is the image of $f'_i$ in $\@O_{B,z}/\mathfrak{n}_z =: \kappa(z)$, the composition
\[
\Spec \kappa(z) \stackrel{z}{\to} B \xrightarrow{pr_1 \circ h''} \#P_k^N
\]
is given by the $(N+1)$-tuple $(\-{f'_0}, \-{f'_1}, \cdots, \-{f'_N})$. of elements in $\kappa(z)$.  Note that $\-{f'_{i_z}} =1$. 

Recall (see equation (\ref{equation: good approximation 1})) that $g_i = f_i( 1 + (p/f_i)(\mu_i + \lambda_i t - \beta_i))$. Let $g'_i = f'_i ( 1 + (p/f_i)(\mu_i + \lambda_i t - \beta_i))$ for every $i$. Thus, we have $g_i = f_{i_z} g'_i$ for every $i$. 

Let $\-{g'_i}$ denote the image of $g'_i$ in $\kappa(z)$. As we observed above, $p/f_i$ is in $\mathfrak{n}_z$. Thus, we see that $\-{f'_i} = \-{g'_i}$ for every $i$. 

The composition 
\[
\Spec \kappa(z) \stackrel{z}{\to} B \xrightarrow{{pr_1 \circ \~h'' \circ \theta'}} \#P_k^N
\]
is given by the $(N+1)$-tuple $(\-{g'_0}, \-{g'_1}, \cdots, \-{g'_N})$. Since $\-{g'_i} = \-{f'_i}$ for every $i$, we see that $pr_1 \circ h''(z) = pr_1 \circ \~h'' \circ \theta'(z)$. 

Similarly, using equation (\ref{equation: good approximation 2}), we can show that $pr_2 \circ h''(z) = pr_2 \circ \~h'' \circ \theta'(z)$ for any point $z \in \tau^{-1}(\@Z(IJ))$. Thus, we conclude that $h''(z) = \~h'' \circ \theta'(z)$ for any $z \in \tau^{-1}(\@Z(IJ))$. 

Now, suppose $z_1$ and $z_2$ are two distinct points of $\~B$ such that $\~\tau(z_1) = \~\tau(z_2) = z$. Then, as $\~\tau$ is an isomorphism on the complement of $\@Z(\~I \~J)$, we see that $z \in \@Z(\~I \~J)$. Recall that $\theta$ maps $\theta^{-1}(\@Z(\~I \~J)) = \@Z(IJ)$ isomorphically onto $\@Z(\~I \~J)$ and that $\theta'$ maps 
$\tau^{-1}(\@Z(I J))$ maps isomorphically onto $\~\tau^{-1}(\@Z(\~I \~J))$. Thus, there exist unique points $y_1$ and $y_2$ in $\psi^{-1}(\@Z(I J))$ such that $\theta'(y_i) = z_i$ for $i = 1,2$. Also, $\~\tau \circ \theta' = \theta \circ \tau$, and since $\theta$ is an injective on $\@Z(IJ)$, we see that $\tau(y_1) = \tau(y_2)$. Since $\chi \circ h'' = h \circ \tau$, we see that $\chi \circ h''(y_1) = \chi \circ h''(y_2)$. Thus, 
\begin{align*}
\chi \circ \~h''(z_1) & = \chi \circ \~h'' \circ \theta'(y_1) \\
                          & = \chi \circ h'' (y_1) \\
                          & = \chi \circ h'' (y_2) \\
                          & = \chi \circ \~h'' \circ \theta'(y_2) = \chi \circ \~h''(z_2) \text{.}
\end{align*}
 
 Note that $\~\tau$ is a proper, birational morphism. Also, as $\Spec R[t]$ is normal, we have $\tau_*(\@O_{\~B}) \cong \@O_{\Spec R[t]}$. Thus, we may apply \cite[Lemma 8.11.1]{EGA2} to conclude that there exists a morphism $\~h: \Spec R[t] \to X$ such that $\chi \circ \~h'' = \~h \circ \~\tau$. 
 
 Now, we will show that $\~h$ can be extended to a morphism $H: \#P^1_R \to X$. (Recall that $\#P^1_R = \Proj R[T_0, T_1]$ contains $\Spec R[t]$ as the open subscheme $\#P^1_R \backslash (\Spec R \times \{\mathbf{\infty}\})$, via the identification $t = T_0/T_1$.) 
 
For $0 \leq i \leq N$, let $G_i(T_0, T_1) \in R[T_0, T_1]$ be defined by
\[
G_i(T_0, T_1) = T_1^{\deg_t(p) + 1} \cdot g_i(T_0/T_1) \text{.}
\]
Thus, each $G_i$ is a homogeneous polynomial of degree $\deg_t(p) + 1$ in $R[T_0, T_1]$. The coefficient of $T_0^{\deg_t(p) + 1}$ in $G_i(T_0, T_1)$ is $\lambda_i$. Recall that there exists an index $i_0$ such that $\lambda_{i_0}$ is a non-zero element of $k$. Thus, $G_{i_0}$ has no zero in the closed subscheme $\@Z(T_1)$ of  $\#P^1_R$. Thus, the rational map $H': \#P^1_R \dashrightarrow \#P^N_k$ defined by the $(N+1)$-tuple $(G_0, \ldots, G_N)$ is defined on an open subscheme of $\#P^1_R$ containing the closed subscheme $\@Z(T_1)$. The restriction of $H'$ to the open subscheme $\#P^1_R \backslash \@Z(T_1) = \Spec R[t]$ is given by the $(N+1)$-tuple $(g_0(t), \ldots, g_N(t))$. Thus, it is the same as the rational map $\~h'$. Thus, $H'$ is defined on the open subscheme $\#P^1_R \backslash \@Z(\~I)$. (Note that the morphism $\@Z(\~I) \hookrightarrow \Spec R[t] \hookrightarrow \#P^1_R$ is closed since $\@Z(\~I)$ is finite over $\Spec R$. The same is true for the closed subschemes $\@Z(\~J)$ and $\@Z(\~I \~J)$ of $\Spec R[t]$. So, we view $\@Z(\~I)$, $\@Z(\~J)$ and $\@Z(\~I \~J)$ as closed subschemes of $\#P^1_R$.)

On the open subscheme $\Spec R[t] \backslash \@Z(\~I\~J)$ of $\Spec R[t]$, where $\psi \circ \~h'$ is well-defined, it agrees with the restriction of the morphism $\~h$. Thus, we see that $\psi \circ H'$ agrees with $\~h$ on $\#P^1_R \backslash (\#Z(\~I\~J) \cup \@Z(T_1))$. Since $\#P^1_R \backslash \@Z(T_1)$ (where $\~h$ is defined) and $\#P^1_R \backslash \@Z(\~I\~J)$ (where $\psi \circ H'$ is defined) form a Zariski open cover of $\#P^1_R$, we see that there exists a morphism $H: \#P^1_R \to X$ extending $\~h$. Now we compute the morphisms $H \circ \sigma_0$ and $H \circ \sigma_{\infty}$ from $\Spec R$ to $X$. 

The morphism $H \circ \sigma_0$ is the same as $\~h \circ \sigma_0$. So we will now compute $h \circ \sigma_0$. We will prove that it is the same as the canonical morphism $\omega_x$. Note that $h \circ \^\sigma_0 = \omega_x$. (Recall from subsection \ref{subsection infinitesimal homotopies} that $\^\sigma_0$ is the morphism $\Spec R \to \Spec R\{t\}$ corresponding to the quotient homomorphism $R\{t\} \to R\{t\}/tR\{t\} \cong R$.) Thus, it is enough to prove that $\~h \circ \sigma_0 = h \circ \^\sigma_0$. For this, it will suffice to show that if $\eta$ denotes the generic point of $\Spec R$, then the two compositions
\[
\xymatrix{
\Spec \kappa(\eta) \ar[r]& {\Spec R} \ar@<-0.5ex>[rr]_{\~h \circ \sigma_0} \ar@<0.5ex>[rr]^{h \circ \^\sigma_0} && X
}
\]
are equal. 

The rational map $h': \Spec R \{t\} \dashrightarrow \#P^N_k$ is defined on the open subscheme $\Spec R\{t\} \backslash \@V(\~I)$ of $\Spec R\{t\}$ and it is represented by the $(N+1)$-tuple $(f_0, \ldots, f_N)$. Since $t$ divides $p$, the image of $f_i$ under the quotient homomorphism $R\{t\} \to R\{t\}/tR\{t\} \cong R$ is $\alpha_i(0)$, i.e. the constant term in the polynomial $\alpha_i(t) \in R[t]$. Since the $f_i$ were chosen to be coprime, we see that at least one of the elements $\alpha_0(0), \ldots, \alpha_N(0)$ is non-zero. Thus, the morphism $h' \circ \^\sigma \circ \eta: \Spec \kappa(\eta) \to \#P^N$ is represented by the $(N+1)$-tuple $(\alpha_0(0), \ldots, \alpha_N(0))$. Similarly, the morphism $\~h' \circ \sigma_0 \circ \eta: \Spec \kappa(\eta) \to \#P^N$ is also represented by the same $(N+1)$-tuple. Composing with $\psi$, we obtain the desired conclusion that the morphism $h \circ \^\sigma_0 \circ \eta$ is equal to $\~h \circ \sigma_0 \circ \eta$. Thus, it follows that $\~h \circ \sigma_0$ is equal to $\omega_x$. 

The restriction of $H'$ to the closed subscheme $\@Z(T_1)$ of $\#P^1_R$ maps $\@Z(T_1)$ to the point $(\lambda_0: \ldots: \lambda_N)$ of $\#P^N_k$. The morphism $H$ agrees with $H'$ on an open subscheme of $\#P^1_R$ containing $\@Z(T_1)$. Thus $H \circ \sigma_{\infty}$ maps $\Spec R$ to the $k$-valued point $\psi((\lambda_0: \ldots: \lambda_N))$ of $X$. This completes the proof.  
 \end{proof}

\end{document}